\documentstyle[12pt]{article}
\newtheorem{Th}{Theorem}
\newtheorem{Prop}[Th]{Proposition}
\newtheorem{Def}[Th]{Definition}
\newtheorem{Rem}[Th]{\rm Remark}
\newtheorem{Cor}[Th]{Corollary}
\newtheorem{Lemma}[Th]{Lemma}
\newcommand{\Proof}{\noindent {\it Proof. }}
\newcommand{\be}{\begin{eqnarray*}}
\newcommand{\ee}{\end{eqnarray*}}
\newcommand{\pkEF}{{\cal P}(^k\!E,F)}

\newcommand{\fin}{\hfill $\Box$\vspace{\baselineskip}}
\newcommand{\series}{\sum_{i=1}^{\infty} x_{i}}

\newcommand{\R}{\bf R}
\newcommand{\C}{\bf C}
\newcommand{\N}{\bf N}
\newcommand{\lra}{\longrightarrow}

\begin{document}

\title{Unconditionally converging polynomials
on Banach spaces\thanks{1991 Mathematics Subject Classification.
Primary: 46B20, 46G20} }

\author{Manuel Gonz\'alez\thanks{Supported in part
by DGICYT Grant PB 91-0307 (Spain)} \\
Departamento de Matem\'aticas \\
Facultad de Ciencias \\
Universidad de Cantabria \\
39071 Santander (Spain)
\and
Joaqu\'\i n M. Guti\'errez\thanks{Supported in part
by DGICYT Grant PB 90-0044 (Spain)} \\
Departamento de Matem\'atica Aplicada \\
ETS de Ingenieros Industriales \\
Universidad Polit\'ecnica de Madrid \\
Jos\'e Guti\'errez Abascal 2 \\
28006 Madrid (Spain) }

\date{}

\maketitle

\begin{abstract}
We prove that weakly unconditionally Cauchy (w.u.C.) series and
unconditionally converging (u.c.) series are preserved under the
action of polynomials or holomorphic functions on
Banach spaces, with natural restrictions in the latter case. Thus it is
natural to introduce the unconditionally converging polynomials,
defined as polynomials taking w.u.C. series into u.c.\ series, and
analogously, the unconditionally converging holomorphic functions.
We show that most of the classes of polynomials
which have been considered in the literature consist of unconditionally
converging polynomials. Then we study several ``polynomial
properties'' of Banach spaces, defined in terms of relations of
inclusion between classes of polynomials, and also some
``holomorphic properties''. We find remarkable differences with the
corresponding ``linear properties''. For example, we show that a Banach
space $E$ has the polynomial property (V) if and only if
the spaces of homogeneous scalar polynomials
${\cal P}(^k\!E)$, $k\in\N$, or the space of scalar holomorphic mappings
of bounded type ${\cal H}_b(E),$ are reflexive. In this case
the dual space $E^*$, like the dual of Tsirelson's space, is reflexive
and contains no copies of $\ell_p$.
\end{abstract}

%\section{Introduction}

In the study of polynomials acting on Banach spaces,
the weak topology is not such a good tool
as in the case of linear operators, due to the bad behaviour
of the polynomials with respect to the weak convergence. For
example,
$$Q : (x_n) \in \ell_2 \lra (x_{n}^{2}) \in \ell_1$$

\noindent
is a continuous polynomial taking a weakly null sequence into a
sequence having no weakly Cauchy subsequences.
In this paper we show that the situation is not so bad for unconditional
series.
Recall that $\sum_{i=1}^{\infty} x_{i}$ is a {\em weakly unconditionally
Cauchy series} (in short a w.u.C. series) in a Banach space $E$ if for
every $f \in E^*$ we have that $\sum_{i=1}^{\infty}
 | f(x_{i}) | < \infty$; and
$\sum_{i=1}^{\infty} x_{i}$ is an {\em unconditionally converging series}
(in short an u.c. series) if every subseries is norm convergent.

We prove that a continuous polynomial takes w.u.C. (u.c.) series into
w.u.C. (u.c.) series. We derive this result from an estimate of the
unconditional norm of the image of a sequence by a homogeneous
polynomial, which is also a fundamental tool in other parts
of the paper, and could be of some interest in itself.

In view of the preservation of w.u.C. (u.c.) series by polynomials,
it is natural to introduce the class ${\cal P}_{uc}$ of
{\em unconditionally converging polynomials} as those taking
w.u.C. series into u.c.\ series. It turns out that most of the
classes of polynomials that have been considered in the literature are
contained in ${\cal P}_{uc}.$ By means of the class of unconditionally
converging polynomials we introduce the {\em polynomial property
$(V)$} for Banach spaces, and we show that
spaces with this property share some of the properties of
Tsirelson's space $T^*.$ In fact, these spaces are reflexive,
and their dual spaces cannot contain copies of
$\ell_p,$ $1<p<\infty.$ This is a consequence of the following
characterization: A Banach space $E$ has the polynomial property
$(V)$ if and only if the space of scalar polynomials
${\cal P}(^k\!E)$ is reflexive for every positive integer $k.$
We also apply ${\cal P}_{uc}$ to characterize the polynomial
counterpart of other isomorphic properties of Banach spaces,
like the Dieudonn\'e property, the Schur property,
and property $(V^*),$ obtaining
remarkable differences with the corresponding linear (usual)
properties.
This is in contrast with the results of \cite{Ry2},
where it is proved that the polynomial Dunford-Pettis property
coincides with the Dunford-Pettis property.

Throughout the paper, $E$ and $F$ will be real or complex Banach
spaces, $B_E$ the unit ball of $E$, and $E^*$ its dual space.
The scalar field will be always $\R$ or $\C$,  the
real or the complex field,
and we will write $\N$ for the set of all natural numbers.
Moreover, ${\cal P}(E,F)$ will stand for the space of all
(continuous) polynomials from $E$ into $F.$
Any $P\in{\cal P}(E,F)$ can be decomposed as a sum of homogeneous
polynomials:
$P = \sum_{k=0}^nP_k,$ with $P_k\in \pkEF ,$ the space of all
k-homogeneous polynomials from $E$ into $F.$

\section{Unconditionally converging polynomials}

In this section we obtain an estimate for the unconditional norm
of the image of a sequence by a homogeneous polynomial,
and we apply it to prove the
preservation of w.u.C. series and u.c.\ series by homogeneous
polynomials. Then we introduce the class of unconditionally converging
polynomials, and compare it with other classes of polynomials that have
appeared in the literature.

In the proof of the estimate, we will need the generalized Rademacher
functions, denoted by  $s_n(t),$ $n\in\N$, which were
introduced in \cite{AG}. These functions are defined as follows:

Fix $2\leq k\in \N$, and let $\alpha_1 = 1, \alpha_2,\ldots ,\alpha_k$
denote the $k^{th}$ roots of unity.

Let $s_1 : [0,1]\rightarrow \C $ be the step function taking the
value $\alpha_j$ on $((j-1)/k,j/k)$ for $j=1,\ldots ,k.$

Then, assuming that $s_{n-1}$ has been defined, define $s_n$
as follows.
Fix any of the $k^{n-1}$ subintervals $I$ of $[0,1]$ used in the
definition of $s_{n-1}.$
Divide $I$ into $k$ equal intervals $I_1,\ldots ,I_k,$ and  set
$s_n(t) = \alpha_j$ if $t\in I_j.$

The generalized Rademacher functions are orthogonal
\cite[Lemma 1.2]{AG}  in the sense that,
for any choice of integers $i_1,\ldots ,i_k;$ $k\geq 2,$ we have
$$\int_{0}^{1} s_{i_1}(t)...s_{i_k}(t) dt =
\left\{\begin{array}{ll}
1, & \mbox{if $i_1 =\cdots = i_k;$} \\
0, & \mbox{otherwise.}
\end{array}
\right. $$

\begin{Lemma}
\label{ineq}
Let $E$ and $F$ be Banach spaces.
Given  $k \in \N$ there exists a constant $C_k$ such that for every
$P \in \pkEF,$ and $x_1,...,x_n \in E$ we have
$$\sup_{|\epsilon_j | \leq 1} \left\|\sum_{j=1}^n \epsilon_j P x_{j}
\right\| \leq
C_k \sup_{|\nu_j | \leq 1}\left\| P\left(\sum_{j=1}^n \nu_j
x_{j}\right) \right\| .$$

In the complex case we can take $C_k = 1$ for every $k$, and in
the real case, $C_k = (2k)^k/k!$.
\end{Lemma}

\Proof
First we assume that $E$ and $F$ are complex spaces. In this case,
both suprema are attained for some
$|\epsilon_j| = |\nu_j| = 1.$

Given $P\in\pkEF,$ we denote by $\hat{P}$ the associated
symmetric k-linear map.
For any $x_1,...,x_n \in E$ and any complex numbers $\epsilon_j$
with $|\epsilon_j| = 1,$
we can find $f\in F^*,$ $\| f\| = 1,$ such that
$$
\left\|\sum_{j=1}^n \epsilon_j P x_{j}
\right\| =
f\left(\sum_{j=1}^n \epsilon_j P x_{j} \right).
$$
Then, taking complex numbers $\delta_j$ such that
$\delta_j^k = \epsilon_j,$ we obtain
$$
\left\|\sum_{j=1}^n \epsilon_j P x_{j}
\right\| =
f\left(\sum_{j=1}^n P(\delta_j x_{j}) \right )=
$$
$$
\int_{0}^{1} \left(\sum_{j_1,\ldots,j_k=1}^n
s_{j_1}(t)\cdots s_{j_k}(t)
f\circ\hat{P}(\delta_{j_1}x_{j_1},\ldots,
\delta_{j_k}x_{j_k})\right) dt =
$$
$$
\int_{0}^{1} f\circ\hat{P}\left(\sum_{j_1=1}^n
\delta_{j_1}s_{j_1}(t)x_{j_1},\ldots, \sum_{j_k=1}^n
\delta_{j_k}s_{j_k}(t)x_{j_k})\right) dt =
$$
$$
\int_{0}^{1} f\circ P\left(\sum_{j=1}^n
\delta_js_j(t)x_j\right) dt \leq
\sup_{|\nu_j|=1} \left\|P\left(\sum_{j=1}^n
\nu_jx_j\right)\right\|.
$$
In this way the proof for the complex case is finished.

Assume now that $E$ and $F$ are real Banach spaces, and let us denote
by $\bar{E}$ and $\bar{F}$ their respective complexifications.
We can extend the multilinear map $\hat{P}\in L(^k\! E,F)$
associated with the polynomial $P$ to a
multilinear map $\hat{Q}\in L(^k\!\bar{E},\bar{F})$ in a
straightforward way, which in the case $k=2$ is given by
$$
\hat{Q}(x_1+iy_1,x_2+iy_2) = \hat{P}(x_1,x_2) + i\hat{P}(y_1,x_2)
+ i\hat{P}(x_1,y_2) - \hat{P}(y_1,y_2),
$$
and the polynomial $Q\in {\cal P}(^k\!\bar{E},\bar{F})$
associated with $\hat{Q}$ is an extension of $P.$ We have
$$
\left\|\sum_{j=1}^n \epsilon_j Px_j\right\| \leq
\sup_{|\nu_j|=1} \left\|Q\left(\sum_{j=1}^n
\nu_jx_j\right)\right\|;
$$
and for complex numbers $\nu_j = a_j + ib_j$ with $|\nu_j|=1,$
we obtain
$$
\left\|Q\left(\sum_{j=1}^n \nu_jx_j\right)\right\|
=
\left\|\sum_{m=0}^k \left( \begin{array}{c} k \\ m \end{array}\right)
i^m\hat{Q}\left(\sum_{j=1}^n b_jx_j\right)^m
\left(\sum_{j=1}^n a_jx_j\right)^{k-m} \right\| =
$$
$$
\left\|\sum_{m=0}^k \left( \begin{array}{c} k \\ m \end{array}\right)
i^m\hat{P}\left(\sum_{j=1}^n b_jx_j\right)^m
\left(\sum_{j=1}^n a_jx_j\right)^{k-m} \right\|.
$$
Moreover, using the polarization formula \cite[Theorem 1.10]{Mu},
we obtain
$$
\left\| \hat{P}\left(\sum_{j=1}^n b_jx_j\right)^m
\left(\sum_{j=1}^n a_jx_j\right)^{k-m} \right\| =
$$
$$
\frac{1}{k!2^k}\left\|\sum_{\epsilon_j = \pm 1}
\epsilon_1\cdots\epsilon_k
P\left((\epsilon_1+\cdots
+\epsilon_m) \left(\sum_{j=1}^n b_jx_j\right) +
(\epsilon_{m+1}+\cdots +\epsilon_k)
\left(\sum_{j=1}^n a_jx_j\right)\right) \right\| =
$$
$$
\frac{1}{k!2^k}\left\|\sum_{\epsilon_j = \pm 1}
\epsilon_1\cdots\epsilon_k k^k
P\left(\sum_{j=1}^n c_j^{\epsilon}x_j\right) \right\|
\leq
\frac{k^k}{k!}\sup_{|c_j|\leq 1}
\left\| P\left(\sum_{j=1}^n c_jx_j\right)\right\|,
$$
where $c_j^{\epsilon} = k^{-1}((\epsilon_1+\cdots +\epsilon_m)bj
+
(\epsilon_{m+1}+\cdots +\epsilon_k)a_j).$
Hence
$$
\left\|\sum_{j=1}^n \epsilon_j Px_j\right\|
\leq
\sum_{m=0}^k \left( \begin{array}{c} k \\ m\end{array}\right)
\frac{k^k}{k!}
\sup_{|c_j| \leq 1} \left\|
P\left(\sum_{j=1}^n c_jx_j\right)\right\|
=
\frac{(2k)^k}{k!} \sup_{|c_j|\leq 1} \left\|P\left(\sum_{j=1}^n
c_jx_j\right)\right\|.
$$
\fin

Next we show that polynomials preserve w.u.C. series
and u.c.\ series.

\begin{Th} \label{series}
Let $E$ and $F$ be Banach spaces and $P \in {\cal P}(E,F).$
Then $P$ takes w.u.C. (u.c.)
series into w.u.C. (u.c.) series.
\end{Th}

\Proof
Recall that a series $\sum_{i=1}^{\infty} x_{i}$ in a
Banach space is w.u.C. if and only if
$\sup_{|\epsilon_i |\leq 1} \|\sum_{i=1}^{\infty}\epsilon_i x_i\|$
is finite; and it is u.c.\ if and only if
$\sup_{|\epsilon_i |\leq 1}\|\sum_{i=n}^{\infty}\epsilon_i x_i \|$
converges to $0$ when $n$ goes to infinity.
Hence, the result is a direct consequence of Lemma \ref{ineq},
since we have
$$
\sup_{|\epsilon_i|\leq 1}\left\|\sum_{i=n}^{\infty}
\epsilon_i Px_i\right\|
\leq C_k \| P\|
\sup_{|\nu_i| \leq 1} \left\|\sum_{i=n}^{\infty}\nu_ix_i\right\|^k.
$$
\fin

Inspired by Theorem \ref{series}, we introduce the following
class of polynomials.

\begin{Def}
{\rm Let $E$ and $F$ be Banach spaces and $k\in\N$.
A polynomial $P\in \pkEF$ is said to be
{\em unconditionally converging}
if it takes w.u.C. series into u.c.\ series.

We shall denote by ${\cal P}_{uc}(^k\!E,F)$ the class of all
unconditionally converging k-homogeneous
polynomials from $E$ to $F$. }
\end{Def}

Note that, in the case $E$ or $F$ contains no copies of $c_0,$
we have that all the w.u.C. series in that space are u.c.; hence
${\cal P}(^k\!E,F) = {\cal P}_{uc}(^k\!E,F).$

The prototype of w.u.C., not u.c.\ series is the unit
vector basis $\{e_n\}$ of $c_0.$ Next lemma characterizes
unconditionally converging polynomials in terms of their
restrictions to subspaces isomorphic to $c_0,$ or the action
on sequences equivalent to $\{e_n\}.$

\begin{Lemma}
\label{uclemma}
Given $P \in \pkEF \setminus {\cal P}_{uc},$ there exists an
isomorphism $i:c_0 \rightarrow E$ such that $\{ (P\circ i)e_n\}$ is
equivalent to $\{e_n\}.$ In particular,
$P\circ i \in{\cal P}(^k\!c_0,F)\setminus {\cal P}_{uc}.$
\end{Lemma}
\Proof
If $P\in\pkEF\setminus {\cal P}_{uc},$ then we can find a w.u.C.
series $\series$ such that
$\sum_{i=1}^{\infty} Px_i$ is not u.c.

Moreover, given a w.u.C., not u.c.\ series $\sum_{i=1}^{\infty}z_i,$
we can construct suitable blocks
$u_k := z_{n_k+1}+\cdots +z_{n_{k+1}}$ such that $\|u_k\|$ is
bounded away from $0.$
Since the series $\sum_{i=1}^{\infty}u_i$ is w.u.C.,
the sequence $(u_k)$ is weakly null;
hence by the Bessaga-Pelczynski selection theorem, it has
a basic subsequence, which is equivalent to the unit vector
basis of $c_0$ \cite[Corollary 5.7]{Di}.

Now, as $\sum_{i=1}^{\infty} Px_{i}$ is w.u.C. but not u.c.,
we can construct a sequence of blocks
$Px_{n_k+1}+\cdots +Px_{n_{k+1}}$ of $(Px_i)$ equivalent to
the unit vector basis of $c_0,$ and it follows from  Lemma~\ref{ineq}
that there exist scalars $c_i$ with $|c_i|\leq 1$ for
every $i\in\N,$ such that the vectors
$P(y_k) = P(c_{n_k+1}x_{n_k+1}+\cdots +c_{n_{k+1}}x_{n_{k+1}})$
are bounded away from $0.$
Since $\sum_{k=1}^{\infty} y_k$ is a w.u.C. series also,
passing to a subsequence if necessary,
we can assume that both sequences $(y_k)$ and $(Py_k)$
are equivalent to $\{ e_k\}.$

The map $i:c_0 \rightarrow E$ defined by  $i(e_n) := y_n$ is
an isomorphism, and
$\sum_{n=1}^{\infty} (P\circ i) e_n$ is not u.c.;
hence $P\circ i \not\in{\cal P}_{uc}(^k\!c_0,F).$
\fin

One of the remarkable properties of the class ${\cal P}_{uc}$
of unconditionally converging polynomials is that it includes
the main classes of polynomials considered in the literature,
as we shall show below.

Recall that a polynomial $P\in\pkEF$ is {\em weakly compact},
denoted by
$P\in {\cal P}_{wco}(^k\!E,F),$ if it takes bounded subsets into
relatively weakly compact subsets.
Moreover, we shall say that $P$ is {\em completely
continuous}, denoted by $P\in {\cal P}_{cc}(^k\!E,F),$ if it
takes weakly Cauchy sequences into norm convergent sequences.
These classes were considered in \cite{Pe} and \cite{Ry2}.

We shall consider also the class ${\cal P}_{cc0}(^k\!E,F)$ of
polynomials which are
{\em completely continuous at 0}, formed by those
$P\in {\cal P}(^k\!E,F)$ taking weakly null sequences
into norm null sequences.
Clearly ${\cal P}_{cc}(^k\!E,F)
\subset {\cal P}_{cc0}(^k\!E,F),$ but in general
(see Proposition~\ref{Schur}) the containment is strict for $k>1$
and $E$ failing the Schur property.

 Recall that $A\subset E$ is said to be a {\it Rosenthal set} if any
sequence $(x_n) \subset A$ has a weakly Cauchy subsequence.
Contrarily to the case of linear operators, a polynomial taking
Rosenthal sets into relatively compact
subsets need not take weakly null sequences into norm null
sequences, as it is shown by the scalar polynomial
$$P : (x_n) \in \ell_2 \lra \sum_{n=1}^\infty x_{n}^{2} \in \R.$$
The converse implication also fails, since for the polynomial
$$Q : (x_n) \in \ell_2 \lra \left(\sum_{k=1}^\infty
\frac{x_{k}}{k}\right) (x_n) \in \ell_2$$
we have that $Q(e_1+e_n) = (1 + 1/n)(e_1+e_n)$ has no convergent
subsequences, although $Q$ takes weakly null sequences into norm
null sequences, because of the factor $(\sum_{k=1}^\infty x_k/k).$

Finally, recall that $A\subset E$ is said to be a {\it Dunford-Pettis
set}~\cite{An} if for any weakly null sequence $(f_n) \subset E^*$
we have
$$\lim_n \sup_{x\in A} |f_n(x)|  = 0.$$

This class of subsets, introduced in \cite{An},
allowed  in \cite{GG1} to define the class ${\cal P}_{wd}$
as follows:

A polynomial $P\in {\cal P}(^k\!E,F)$ belongs to ${\cal P}_{wd}$
if and only if its restriction to any Dunford-Pettis subset of $E,$
endowed with the inherited weak topology, is continuous.

\begin{Prop}
\label{cc0}
A polynomial $P\in\pkEF$ belongs to ${\cal P}_{uc}$ in the following
cases:

(a)  $P\in {\cal P}_{cc0}.$

(b) $P$ takes Rosenthal subsets of $E$ into relatively compact
subsets of $F.$

(c)  $P\in {\cal P}_{wd}.$

(d)  $P\in {\cal P}_{wco}.$
\end{Prop}
\Proof
The result in the cases (a) and (b) is an immediate consequence of
Lemma~\ref{uclemma}, since the unit vector basis of $c_0$ is a
weakly null sequence which forms a non relatively compact set.

Case (c) follows from Lemma~\ref{uclemma} also, since given
$P \in \pkEF \backslash {\cal P}_{uc},$ and an isomorphism
$i:c_0 \rightarrow E$ such that
$P\circ i \not\in  {\cal P}_{uc}(^k\!c_0,F),$ we have that
$\{ ie_n\}$ is a Dunford-Pettis set of $E$ on which $P$ is not
weakly continuous; hence $P \not\in {\cal P}_{wd}.$

Finally, by \cite[Lemma 3.11 and Theorem
3.13]{GG1}, we have ${\cal P}_{wco}(^k\!E,F) \subseteq
{\cal P}_{wd}(^k\!E,F)$; hence (d) follows from (c).
\fin

\section{Polynomial properties}

Banach spaces with the {\it polynomial Dunford-Pettis property}
were introduced in \cite{Pe} as the spaces $E$ such that
weakly compact polynomials from $E$ into any Banach space
are completely continuous; i.e.,${\cal P}_{wco}(^k\!E,F)\subseteq
{\cal P}_{cc}(^k\!E,F)$ for any $k\in\N$ and $F.$
Later Ryan \cite{Ry2} proved that it coincides with the usual
{\em Dunford-Pettis property,} which admits the
same definition in terms of linear operators $(k=1).$
On the other hand, Pelczynski~\cite{Pe2} introduced Banach
spaces with {\em property (V)} as the spaces $E$ such that
unconditionally converging operators from $E$ into any Banach
space are weakly compact.

In this section, by means of the class ${\cal P}_{uc}$ of
unconditionally converging polynomials,
we introduce and study the polynomial property (V)
and other polynomial versions of properties of Banach spaces:
Dieudonn\'e property, Schur property, and property $(V^*).$
We show that, in contrast to the
case of the Dunford-Pettis property, property (V) is very
different from the polynomial property (V), since spaces with this
property are analogous to Tsirelson's space $T^*.$ For the other
polynomial properties, we show that sometimes the polynomial
and the linear properties coincide,
and sometimes not, with a general tendency of the polynomial
property to imply the absence of copies of $\ell_1$ in
the space. Moreover, we obtain additional results relating
${\cal P}_{uc}$ and other classes of polynomials.

\begin{Def}
{\rm A Banach space $E$ has the {\em polynomial property (V)}
if for every $k$ and $F$ we have
${\cal P}_{uc}(^k\!E,F)\subseteq {\cal P}_{wco}(^k\!E,F).$ }
\end{Def}

It was shown in \cite{Pe2} that $C(K)$ spaces enjoy property (V).
Next Lemma
shows that this is not the case for the polynomial property.

\begin{Lemma}
\label{polV}
Given a Banach space $E$, if ${\cal P}_{uc}(^k\!E,E) \subseteq
{\cal P}_{wco}(^k\!E,E)$ for some $k>1,$ then
$E$ contains no copies of $c_0.$

\end{Lemma}

\Proof Assume $E \supset c_0,$ and take a sequence $(x_n) \subset E$
equivalent to the unit vector basis of $c_0.$
We select $f \in E^*$ such that $f(x_1) = 1,$ $f(x_i) = 0$ for
$i>1,$ and define
$$P:x\in E \longrightarrow f(x)^{k-1}x \in E.$$
Note that $P \in {\cal P}_{uc}(^k\!E,E),$ since for every w.u.C.
series $\series$ in $E,$ we have
$$\sum_{i=1}^{\infty} \| P(x_{i})\| \leq \sup_{j\in\N}\| x_j\|
\left( \sum_{i=1}^{\infty} | f(x_{i}) |^{k-1}\right) \leq
\sup_{j\in\N}\| x_j\| \left(\sum_{i=1}^{\infty} |
f(x_{i}) |\right)^{k-1} <
\infty .$$

However,
$P\not\in {\cal P}_{wco},$ since $P(x_1+...+x_n) = x_1+...+x_n,$
and
$(x_1+...+x_n)_{n\in \N}$ is a weakly Cauchy sequence having
no weakly convergent subsequences.
\fin

\begin{Rem}
\label{c0remark}
The proof of Lemma~\ref{polV} gives also the following facts:

(a) Given a Banach space $E$ containing a copy of $c_0,$ for every
$k>1$ there exists $P_k \in {\cal P}_{uc}(^k\!E,E)$ (even
taking w.u.C. series into absolutely converging series) which
does not take Rosenthal sets into relatively weakly compact sets.

(b) If ${\cal P}_{cc0}(^k\!E,E) \subseteq
{\cal P}_{wco}(^k\!E,E)$ for some $k>1,$ then $E$ contains no
copies of $c_0,$ which is in contrast with the linear case
$(k=1)$ too.

\end{Rem}

It is well-known (and can be easily derived from
the case $k=1$ in Proposition~\ref{Pucc0} below) that every
unconditionally converging operator
$T:c_0 \rightarrow F$ is compact. In contrast, Lemma~\ref{polV}
shows that for $k>1$  polynomials
$P\in {\cal P}_{uc}(^k\!c_0,F)$ are not always weakly compact.
However, these polynomials have restrictions to finite codimensional
subspaces with arbitrarily small norm.

\begin{Prop}
\label{Pucc0}
For any $P\in {\cal P}_{uc}(^k\!c_0,F)$ we have
$$\lim_{n\rightarrow\infty } \left\|
P|_{[ e_n,e_{n+1},\ldots ] }\right\| = 0.$$
\end{Prop}
\Proof If there exists $\delta >0$ such that
$\| P|_{[ e_n,e_{n+1},\ldots ]
}\| >\delta$ for all $n\in\N,$ then we can construct blocks
$u_i := a_{n_i +1}e_{n_i +1}+\cdots +a_{n_{i+1}}e_{n_{i+1}},$
with $n_1<n_2<\cdots ,$ such that
$\| u_i\| = 1$ and $\| Pu_i\| >\delta.$

Then $\sum_{i=1}^{\infty}u_i$ is a w.u.C. series, but
$\sum_{i=1}^{\infty}Pu_i$ is not u.c.; hence
$P\not\in {\cal P}_{uc}.$ \fin

 It has been shown \cite{AAD} that a Banach space $E$ such
that ${\cal P}(^k\!E,{\bf K}) \equiv {\cal P}(^k\!E)$
is reflexive for every $k\in\N$ has many of the properties
of Tsirelson's space $T^*$ \cite{Ts}. In fact,
$E$ must be reflexive, and the dual space $E^*$ cannot
contain copies of $\ell_p$ $(1<p<\infty).$ Note also that
${\cal P}(^kT^*)$ is reflexive for every ${k\in\N}$ \cite{AAD}.
Next we present a characterization of the spaces $E$
such that ${\cal P}(^kE)$ is reflexive for some $k>1$
in terms of the class ${\cal P}_{uc}$ of polynomials.

Given $P\in {\cal P}(^k\!E,F),$ we consider the associated conjugate
operator defined by
$$P^* : f\in F^* \longrightarrow f\circ P \in {\cal P}(^k\!E).$$
Moreover, we need the fact that for every Banach space $E,$ the space
$\Delta^k_{\pi}E,$
defined as the closed span of $\{ x\otimes\cdots\otimes x : x\in E \}$
in the projective tensor product
$\hat{\otimes}^k_{\pi}E,$ is a predual of the
space of scalar polynomials ${\cal P}(^k\!E)$ \cite{Ry1}.

\begin{Th}
\label{ThpolV}
Given $k\in\N,$ $k>1,$ and a Banach space $E,$ we have that
${\cal P}_{uc}(^k\!E,F) \subseteq {\cal P}_{wco}(^k\!E,F)$ for any $F$
if and only if
${\cal P}(^k\!E)$ is reflexive.

%Given $k\in\N,$ $k>1,$ and a Banach space $E,$
%the following assertions are equivalent:

%(a)  for any $F,$ we have ${\cal P}_{cc0}(^k\!E,F) \subseteq
%{\cal $P}_{wco}(^k\!E,F);$

%(b)  for any $F,$ we have ${\cal P}_{uc}(^k\!E,F) \subseteq
%{\cal P}_{wco}(^k\!E,F);$

%(c) the space ${\cal P}(^k\!E)$ is reflexive. $(a) \Rightarrow (c)$

\end{Th}

\Proof  Assume ${\cal P}_{uc}(^k\!E,F) \subseteq
{\cal P}_{wco}(^k\!E,F)$ for any $F.$ By Lemma~\ref{polV}
we have that $E$ contains no copies of $c_0.$ Then
$${\cal P}(^k\!E,F) = {\cal P}_{uc}(^k\!E,F) =
{\cal P}_{wco}(^k\!E,F)$$ for any $k$ and $F.$
Since there exists an isomorphism between the space of polynomials
${\cal P}(^k\!E,F)$ and the space of operators
$L(\Delta^k_{\pi}E,F)$ which takes the weakly
compact polynomials onto the weakly compact operators
\cite{Ry1}, we obtain that all operators in
$L(\Delta^k_{\pi}E,F)$ are weakly
compact. Then  $\Delta^k_{\pi}E$ is reflexive;
hence ${\cal P}(^k\!E) \cong (\Delta^k_{\pi}E)^*$ is reflexive.

Conversely, since $P\in {\cal P}_{wco}$ if and only if the operator
$P^*$ is weakly compact \cite[Proposition 2.1]{Ry3},
if ${\cal P}(^k\!E)$ is reflexive, then we have that any
$P\in {\cal P}(^k\!E,F)$ belongs to
${\cal P}_{wco},$ and the result is proved. \fin

\begin{Cor}
A Banach space $E$  has the polynomial property (V)
if and only if ${\cal P}(^k\!E)$ is reflexive for every $k\in\N$.
\end{Cor}

We will need the following well-known characterization of
Banach spaces containing no copies of $\ell_1.$
We include a proof for the sake of completeness.

\begin{Lemma}
\label{l1l2}
A Banach space $E$ contains a copy of $\ell_1$ if and only if there
exists a completely continuous surjection from $E$ onto $\ell_2.$

\end{Lemma}

\Proof Assume $E$ contains a copy of $\ell_1,$ and let $q$
denote a surjective operator from $\ell_1$ onto $\ell_2.$
Since $q$ is absolutely summing \cite[Theorem 2.b.6]{LT},
it factors through a space $L_{\infty}({\mu}),$
which has the extension property and the Dunford-Pettis property.
Then the operator from $\ell_1$ into $L_{\infty}({\mu})$
can be extended to an operator $A$ from
$E$ into $L_{\infty}({\mu}),$ and the operator $B$ from
$L_{\infty}({\mu})$ onto $\ell_2$ is completely continuous;
hence $BA$ is a completely continuous,
surjective operator from $E$ onto $\ell_2.$

Conversely, if $Q$ is a completely continuous, surjective
operator from $E$ onto $\ell_2,$ and we take a bounded
sequence $(x_n)$ in $E$ such that $\{Qx_n\}$ is the
unit vector basis of $\ell_2,$ then $(x_n)$ cannot have a
weakly Cauchy subsequence; hence, by Rosenthal's theorem,
it has a subsequence equivalent to the unit vector
basis of $\ell_1.$ \fin

In relation with the reflexivity of ${\cal P}(^k\!E),$
the problem of when this space contains a copy of $\ell_{\infty}$
has received some attention.
It has been considered in \cite{ACL}, for instance,
in connection with the so called property (RP) of polynomials.
We give an answer that includes the case
$\ell_1 \subseteq E$ (see Lemma~\ref{l1l2}).

\begin{Prop} If $E$ has a quotient isomorphic to $\ell_2,$
then for any integer $k>1$ the space ${\cal P}(^k\!E)$
contains a copy of $\ell_{\infty}.$
\end{Prop}

\Proof For every $a\equiv (a_n)\in \ell_{\infty},$ we consider
the polynomial
$$P_a : (x_n)\in\ell_2 \rightarrow \sum_{i=1}^{\infty} a_i x_i^k.$$

We have $P_a \in {\cal P}(^k\!\ell_2),$ and $\|P_a\|
= \|(a_n)\|_{\infty}.$
Then, the map
$$a\in \ell_{\infty} \rightarrow P_a \in {\cal P}(^k\!\ell_2)$$
defines a linear isometry from $\ell_{\infty}$ into
${\cal P}(^k\!\ell_2).$ Now, if $q : E \rightarrow \ell_2$
is a quotient map, we have that the map
$$a\in \ell_{\infty} \rightarrow P_a\circ q \in {\cal P}(^k\!E)$$
is an isomorphism from $\ell_{\infty}$ into ${\cal P}(^k\!E).$ \fin

Extending the definition for operators,
we shall say that a polynomial
$P\in {\cal P}(^k\!E,F)$ is {\it weakly completely
continuous}, denoted by $P \in {\cal P}_{wcc}(^k\!E,F),$ if it
takes weakly Cauchy sequences into weakly convergent sequences.

A Banach space $E$ has the {\em Dieudonn\'e property} if weakly
completely continuous operators from $E$ into
any Banach space are weakly compact.
Grothendieck \cite[3.1]{Gr} introduced this property and
proved that $C(K)$ spaces enjoy it. Next
result shows that the polynomial Dieudonn\'e property is
equivalent to the
absence of copies of $\ell_1$ in the space.

\begin{Prop}
\label{nl1}
For a Banach space $E$ the following properties are equivalent:

(a)  $E$ contains no copies of $\ell_1.$

(b) ${\cal P}_{wcc}(^k\!E,F)\subseteq {\cal P}_{wco}(^k\!E,F)$ for
any $k$ and $F.$

(c) ${\cal P}_{cc}(^k\!E,F)\subseteq {\cal P}_{wco}(^k\!E,F)$ for
any $k$ and $F.$

(d) ${\cal P}_{cc}(^k\!E,F)\subseteq {\cal P}_{wco}(^k\!E,F)$ for
some nonreflexive $F$ and some $k > 1.$

\end{Prop}

\Proof $(a) \Rightarrow (b)$ Assume $E$ contains no copies of
$\ell_1,$ and let $P \in {\cal P}_{wcc}(^k\!E,F).$
Since any bounded sequence $(x_n) \subset E$ has a weakly
Cauchy subsequence, we have that $(Px_n)$ has a weakly convergent
subsequence; hence $P \in {\cal P}_{wco}.$

$(b) \Rightarrow (c) \Rightarrow (d)$  are trivial.

$(d) \Rightarrow (a)$ Assume $E$ contains a copy of $\ell_1,$
and $F$ is nonreflexive.
We take a sequence $(y_n) \subset B_F$ having no weakly
convergent subsequences. By Lemma~\ref{l1l2},
 we can also take a
completely continuous surjection $T : E \rightarrow \ell_2.$
Now, if $k>1,$  $Q$ is the polynomial from $\ell_2$ into
$\ell_1$ defined by $Q(x_i) := (x_i^k),$ and
$S:\ell_1\rightarrow F$ is the operator defined by $Se_n := y_n,$
where ${e_n}$ is the unit vector basis of $\ell_1,$
then $S\circ Q\circ T \in {\cal P}_{cc}(^k\!E,F)$, but
$S\circ Q\circ T \not\in {\cal P}_{wco},$ because there
exists a bounded sequence
$(x_n)\subset E$ such that $S\circ Q\circ T x_n = y_n.$  \fin

\begin{Cor}
${\cal P}_{wcc}(^k\!E,F) \subseteq {\cal P}_{uc}(^k\!E,F)$
for any $k\in\N$.
\end{Cor}
\Proof If $P \in {\cal P}(^k\!E,F) \backslash {\cal P}_{uc},$ by
Lemma~\ref{uclemma}, there exists an isomorphism
$i:c_0 \rightarrow E$ such
that $P\circ i \in {\cal P}(^k\!c_0,F) \backslash
{\cal P}_{uc}.$ Then, by Proposition~\ref{cc0},
$P\circ i \not\in {\cal P}_{wco}(^k\!c_0,F),$ and by
Proposition~\ref{nl1},  $P\circ i\not\in {\cal P}_{wcc};$
hence $P\not\in {\cal P}_{wcc}.$ \fin

\begin{Rem}
It follows from Proposition~\ref{nl1} that, for any $k>1,$ there is a
polynomial $P\in {\cal P}_{cc}(^k\!\ell_{\infty},c_0)$ which is not
weakly compact. However, any operator from $\ell_{\infty}$
into $c_0$ is weakly compact and thereby completely
continuous, since $\ell_{\infty}$ has the Dunford-Pettis property.

Then the question arises whether every polynomial from $\ell_{\infty}$
into $c_0$ is completely continuous.
\end{Rem}

As a complement of Theorem~\ref{ThpolV} we have the following

\begin{Th}
\label{ThpolV2}
Given $k\in\N,$ $k>1,$ and a Banach space $E,$ we have that
${\cal P}_{cc0}(^k\!E,F) \subseteq {\cal P}_{wco}(^k\!E,F)$
for any $F$ if and only if
${\cal P}(^{k-1}\!E)$ is reflexive.
\end{Th}

\Proof First we assume that ${\cal P}_{cc0}(^k\!E,F)
\subseteq {\cal P}_{wco}(^k\!E,F),$ and
as in the proof of Theorem~\ref{ThpolV}, it is enough to prove that
${\cal P}(^{k-1}\!E,F) = {\cal P}_{wco}(^{k-1}\!E,F).$

By Proposition~\ref{nl1} we have that $E$ contains no copies
of $\ell_1$. Then, if there exists
$P\in {\cal P}(^{k-1}\!E,F),$ $P\notin {\cal P}_{wco},$
we can find a weakly Cauchy sequence $(x_n)\subset E$ such
that $(Px_n)$ has no weakly convergent subsequence;
then, it is not relatively weakly compact.
Since the class of relatively weakly compact sequences is
closed in the space of bounded sequences,
we can assume also that $(x_n)$ is not weakly null.

Taking $\phi\in E^*$ such that
$\phi (x_n) \rightarrow\lambda \neq 0,$
we define a polynomial
$Q\in {\cal P}(^k\!E,F)$ by $Q(x) := \phi (x) P(x).$ Since
$Q\in {\cal P}_{cc0},$ we have that $Q$ is weakly compact.
Then there exists a subsequence $(y_k)$ of $(x_n)$ such that
$Q(y_k)$ is weakly convergent to $z\in E$; therefore,
$P(y_k)$ is weakly convergent to $\lambda^{-1}z.$ Contradiction.

Conversely, if ${\cal P}(^{k-1}\!E)$ is reflexive, we have that
${\cal P}(^l\!E,F) = {\cal P}_{wco}(^l\!E,F)$ for all $l<k$.
Moreover, since $E$ is reflexive,
given a bounded sequence $(x_n)\subset E$,
we can assume, passing to a subsequence, that it is weakly
convergent to some $x\in E.$ Then, given
$P\in {\cal P}_{cc0}(^k\!E,F),$ using the associated
multilinear map $\hat{P}$ we write
$$
Px_n = \hat{P}(x_n-x+x,...,x_n-x+x)
=
\sum_{l=1}^{k-1} \hat{P}(x_n-x)^l(x)^{k-l} + P(x_n-x) + P(x).
$$
Since for $l<k$ the polynomials $Q_l \in {\cal P}(^l\!E,F)$
defined by
$Q_l(y) := \hat{P}(y)^l(x)^{k-l}$ are weakly compact,
and $P(x_n-x)$ converges to $0$,we obtain that
$P(x_n)$ has a weakly convergent subsequence; hence
$P\in {\cal P}_{wco}.$ \fin

\begin{Rem}
In order to compare Theorems \ref{ThpolV} and \ref{ThpolV2}, we
observe that for the sequence spaces $\ell_p$ the space of
polynomials ${\cal P}(^k\!\ell_p)$
is reflexive if and only if $k<p<\infty$.

In fact, it was proved in \cite[Corollary 4.3]{Pe3} that
for $k<p$, all polynomials in
${\cal P}(^k\!\ell_p)$ are completely continuous;
hence, using a result of
\cite{Ry1} (see \cite[Proposition 3]{AAD}), we conclude that
${\cal P}(^k\!\ell_p)$ is reflexive. For $1<p\leq k$ it
is not difficult to show that
${\cal P}(^k\!\ell_p)$ contains a copy of $\ell_{\infty}.$

\end{Rem}

Recall that a Banach space $E$ has the {\em Schur property}
if weakly convergent sequences in $E$ are norm convergent;
equivalently, weakly Cauchy sequences are norm convergent.
It is an immediate consequence of the
definition that $E$ has the Schur property if and only if
${\cal P}(^k\!E,F) = {\cal P}_{cc}(^k\!E,F)$ for any $k$ and $F.$
Next we give
some other polynomial characterizations of Schur property.

\begin{Prop}
\label{Schur}
For a Banach space $E$ the following properties are equivalent:

(a)  $E$ has the Schur property.

(b) ${\cal P}_{uc}(^k\!E,F)\subseteq {\cal P}_{cc}(^k\!E,F)$ for
any $k$ and $F.$

(b') ${\cal P}_{cc0}(^k\!E,F)\subseteq {\cal P}_{cc}(^k\!E,F)$ for
any $k$ and $F.$

(c) ${\cal P}_{uc}(^k\!E,E)\subseteq {\cal P}_{cc}(^k\!E,E)$
for some $k>1.$

(c') ${\cal P}_{cc0}(^k\!E,E)\subseteq {\cal P}_{cc}(^k\!E,E)$
for some $k>1.$

\end{Prop}

\Proof (a) $\Rightarrow$ (b) is immediate.

(b) $\Rightarrow$ (c) $\Rightarrow$ (c') and (b) $\Rightarrow$ (b')
$\Rightarrow$ (c') follow
from ${\cal P}_{cc0}\subseteq{\cal P}_{uc}$ (see
Proposition~\ref{cc0}).

(c') $\Rightarrow$ (a) Assume $E$ fails the Schur property.
We take $x_0\in E,$ with $\| x_0\| = 1,$ and
$f\in E^*$ such that $f(x_0) = 1.$ Since the kernel
of $f$ fails the Schur property also, there exists a weakly
null, normalized sequence $(x_n) \subset ker(f).$
Now, for every $k>1$ we can define a
polynomial $P\in {\cal P}(^k\!E,E)$ by
$$P : x\in E \longrightarrow f(x)^{k-1}x\in E.$$
We have $P\in {\cal P}_{cc0},$
and $x_0+x_n \stackrel{w}{\rightarrow} x_0,$
but $P(x_0+x_n) = x_0+x_n$ does not converge in norm
to $Px_0;$ hence $P \not\in {\cal P}_{cc}.$   \fin

Recall that a Banach space is said to have the
{\em hereditary Dunford-Pettis property}
if any of its subspaces has the Dunford-Pettis property.

\begin{Prop}
\label{hDPP}
If a Banach space $E$ has the hereditary Dunford-Pettis property,
then ${\cal P}_{uc}(^k\!E,F) \subseteq
{\cal P}_{cc0}(^k\!E,F)$ for any $k$ and $F.$
\end{Prop}
\Proof Given $P\in {\cal P}_{uc}(^k\!E,F),$ since $E$ has
the hereditary Dunford-Pettis property, every normalized weakly
null sequence in $E$ has a subsequence equivalent to the unit
vector basis of $c_0$ \cite[Proposition 2]{Ce}, which is taken
into a norm null sequence by $P.$ Thus, every weakly null sequence
$(x_n) \subset E$ has a subsequence $(x_{n_i})$ such that
$(Px_{n_i})$ is norm null; hence $P\in {\cal P}_{cc0}.$ \fin

\begin{Rem}
We do not know if the converse of Proposition~\ref{hDPP} is true.
\end{Rem}

Another property of Banach spaces defined in terms of series is the
property $(V^*),$ introduced by Pelczynski in \cite{Pe2}.
Recall that a subset $A \subset E$ is said to be a $(V^*)$ $set$
if for every w.u.C. series
$\sum_{n=1}^{\infty}f_n$ in $E^*$ we have
$$\lim_n \sup_{x\in A}  | f_n(x) | = 0.$$
A Banach space $E$ has  $property$ $(V^*)$
if every $(V^*)$ set in $E$  is relatively weakly
 compact; equivalently,
if any operator $T\in L(F,E),$ with
unconditionally converging conjugate $T^*$ is weakly compact.

Next we shall show that the polynomial version of the last
formulation coincides with property $(V^*).$
We shall denote by ${\cal P}_{uc*}(^k\!F,E)$ the class of all
polynomials $P\in{\cal P}(^k\!F,E)$ such that $P^*$ is
unconditionally converging.

\begin{Prop}
\label{P*uc}
Given $P\in{\cal P}(^k\!F,E),$ we have that $P^*$ is unconditionally
converging if and only if $P(B_F)$ is a $(V^*)$ set.
\end{Prop}

\Proof Assume $P^*$ is unconditionally converging and
$\sum_{n=1}^{\infty}f_n$ is a w.u.C. series in $E^*.$
We have that $\sum_{n=1}^{\infty}P^*f_n$ is an
u.c.\ series; in particular, $\|P^*f_n\| \rightarrow 0.$ Then
$$\lim_n \sup_{x\in PB_F} | f_n(x) | =
\lim_n \sup_{y\in B_F} | (P^*f_n)y | = 0;$$
hence,
$PB_F$ is a $(V^*)$ set.

Conversely, if $PB_F$ is a $(V^*)$ set, it follows in an
analogous way that $\|P^*f_n\| \rightarrow 0$ for every w.u.C.
series $\sum_{n=1}^{\infty}f_n$ in $E^*;$ and using
Lemma~\ref{uclemma} in the case $k=1,$ we conclude that
$P^*$ is unconditionally converging. \fin

\begin{Prop}
For a Banach space $E$ the following properties are equivalent:

(a) $E$ has property $(V^*).$

(b) For any $k$ and any $F,$ we have ${\cal P}_{uc*}(^k\!F,E)
\subseteq {\cal P}_{wco}(^k\!F,E).$

(c) For some $k,$ we have  ${\cal P}_{uc*}(^k\!\ell_1,E)
\subseteq {\cal P}_{wco}(^k\!\ell_1,E).$

\end{Prop}

\Proof $(a) \Rightarrow (b)$ Assume $E$ has property $(V^*)$
and let $P\in{\cal P}_{uc*}(^k\!F,E).$ Then $PB_F$ is a
$(V^*)$ set; hence it is relatively weakly compact, and we
conclude  $P\in {\cal P}_{wco}.$

$(b) \Rightarrow (c)$ is trivial.

$(c) \Rightarrow (a)$ Assume $E$ fails property $(V^*).$
Then there exists a bounded sequence $(x_n)\subset E,$
having no weakly convergent subsequences,
such that $\{ x_n\}$ is a $(V^*)$ set.
Now, for any $k,$ we define $P\in {\cal P}(^k\!\ell_1,E)$ as
follows:
$$P : (t_i) \in \ell_1 \longrightarrow
\sum_{i=1}^{\infty}t_i^k x_i.$$

Since $PB_{\ell_1}$ is contained in the absolutely convex,
closed hull of $\{ x_n\}$ it is a $(V^*)$ set; hence, by
Proposition~\ref{P*uc}, $P\in {\cal P}_{uc*}.$
However, as $Pe_i= x_i$ for every $i\in\N ,$
where $\{e_i\}$ stands for the unit vector basis of
$\ell_1,$ we have that $P\not\in {\cal P}_{wco}.$ \fin

\begin{Rem}

Using the Taylor expansion, it is possible to show that
holomorphic mappings preserve (locally)
w.u.C. series and u.c.\ series, and a
holomorphic map $f:E\rightarrow F$ is unconditionally converging
if and only if $f(0)=0$ and the homogeneous polynomials given
by the derivatives of $f$ at the origin,
$d^kf(0)$ ($k\in\N),$ are unconditionally converging.
Hence, it follows that the "holomorphic"
property (V) coincides with the polynomial
property (V), obtaining in this way a
characterization of Banach spaces $E$
such that the space ${\cal H}_b(E)$ of
holomorphic mappings of bounded type on $E$ is reflexive.

\end{Rem}

\noindent
{\bf Acknowledgement.} The authors are indebted to Professor
J. Diestel for suggesting the study of the polynomial property (V).

\end{document}